\documentclass[11pt]{article}
         
        \usepackage[dvipsnames]{xcolor}
        \usepackage{keyval}
        \usepackage{diagram}
        \usepackage{geometry}                
        \usepackage{amssymb}
        \usepackage{amsmath}
        \usepackage{enumerate}
        \usepackage{subcaption}
        \usepackage{sectsty}
        \usepackage{float}
        \usepackage{graphicx}
        \usepackage{subfig}
\usepackage{bm}
        
\usepackage[LSBC3,T1]{fontenc}
\usepackage{chessboard}

        \newcommand{\qed}{\mbox{$\Box$}\vspace{\baselineskip}}
         
        \newtheorem{theorem}{Theorem}[section]

        \newtheorem{corollary}[theorem]{Corollary}
        \newtheorem{ques}[theorem]{Question}
        \newtheorem{definition}[theorem]{Definition}
        
        \newtheorem{claim}[theorem]{Claim}

        \title{Exploring mod 2 $n$-queens games}
\author{Tricia Muldoon Brown\\
  \textit{Armstrong State University}
  \and
  Abrahim Ladha\\
  \textit{Armstrong State University}
}
\date{}

\begin{document}
\maketitle

\begin{abstract}
We introduce a two player game on an $n\times n$ chessboard where queens are placed by alternating turns on a chessboard square whose availability is determined by the number of queens already on the board which can attack that square modulo two.  The game is explored along with some variations and its complexity.
\end{abstract}
         
        \setboolean{piececounter}{false}
        \setboolean{showcomputer}{false}
         
        \section{Introduction}
         
        The $n$-queens problem is the problem of placing $n$ queens on an $n$ by $n$ chessboard such that none are attacking each other.  This puzzle originated as the $8$-queens problem played on a standard $8\times 8$ chessboard and was proposed by Max Bezzel in 1848~\cite{Bezzel}.  The $8$-queens problem and the more general $n$-queens problem attracted the interest of notable mathematicians of that time.  In 1850, Nauck was the first to publish all 92 solutions for the standard $8\times 8$ chessboard~\cite{Nauck} and in 1874 Pauls gives the first set of solutions for the general $n$-queens problem published in two articles~\cite{Pauls, Pauls_2}.
             
The problem of finding more solutions in various dimensions has continued to hold the interest of mathematicians and computer scientists.  Solutions sets have been found using graph theory, magic squares, Latin squares, and group theory, among other techniques.  In a brute force approach, solution sets for a given $n$ can be found use the backtracking algorithm, a standard technique taught in computer science classes.  In more modern times, a player can find numerous mobile apps to test his ability to find solutions on $8\times 8$ or other dimensional boards.  

In this paper, we look at a game variant of the $n$-queens problem that can be played on chessboards of varying dimensions.  The basic version of the $n$-queens game is described by Noon in~\cite{Noon}, who, observing that not all placements of queens will lead to a full solutions with $n$ queens, suggests a two player game where each player successively places queens in non-attacking positions.  The first player who cannot place a queen loses the game.  In the next section, we introduce a modification of this $n$-queens game.  Sections~\ref{Complete} and~\ref{Locked} discuss two states a game board may take.  Then Section~\ref{Alternate} will discuss alternate versions of the game, while Section~\ref{Graph} looks at the complexity of the game.  
         
\section{The mod 2 \textit{n}-queens game}
We begin by observing that counting all possible arrangement of queens which satisfy Noon's $n$-queens game can be quite complicated.  This can be illustrated by using a gradient to indicate the number of queens attacking each square on the chessboard, with white indicating an open square and darker shades indicating an increasing number of queens attacking that square.  Figure~\ref{gradient_board} displays an example of such a chessboard.  
\begin{figure}
\begin{subfigure}[t]{.46\textwidth}
\centering
 \definecolor{mybgcolor}{RGB}{255,255,255}
    \definecolor{mygridcolor}{RGB}{0,0,0}
    \definecolor{myhighlightcolor}{RGB}{102,102,102}
    \definecolor{gray1}{RGB}{204,204,204}
    \definecolor{gray2}{RGB}{153,153,153}
\setchessboard{boardfontencoding=LSBC3,
vlabel=false,
hlabel=false,
showmover = false,
    maxfield = h8,
    boardfontsize=30pt,
    boardfontfamily=skaknew,
    pgfstyle=border,
    color=mygridcolor,
    linewidth=0.5pt,
    markboard,
    pgfstyle=color,
    color=mybgcolor,
    backboard,
    color=myhighlightcolor,
    backregions={c8-d8,c6-c6,d5-d5,f5-f5,a4-a4,d2-d2,g2-g2},
    color=gray1,
    backregions={h1-h6,h8-h8,f6-f8,d6-e6,a8-b8,a6-a6,e5-e5,c4-d4,f4-f4,h4-h4,a3-c3,e3-g3,b2-c2,e2-f2,a1-b1,g1-g1,d7-d7,g5-g5},
    color=gray2,
    backregions={a7-c7,e7-e8,g6-g8,h7-h7,g4-g4,c1-d1,a5-c5,e4-e4,d3-d3,a2-a2},
   addblack={qa8,qd7,qc2,qg5},
}
    \makeatletter
    \let\color@endgroupORI\color@endgroup
    \def\color@endgroup{\color@endgroupORI\pgfsetfillopacity{1}}
    \def\cfss@whitefieldmaskcolor{\pgfsetfillopacity{0}\color{white}}
    \def\cfss@blackfieldmaskcolor{\pgfsetfillopacity{0}\color{black}}
    \def\cfss@whitefieldcolor{\pgfsetfillopacity{0}\color{white}}
    \def\cfss@blackfieldcolor{\pgfsetfillopacity{0}\color{black}}
    \makeatother
    \scalebox{.5}{
    \chessboard}
\caption{A chessboard with a gradient indicating attacking queens.}
\label{gradient_board}
\end{subfigure}
\qquad
\begin{subfigure}[t]{.46\textwidth}
\centering
\definecolor{mybgcolor}{RGB}{255,255,255}
    \definecolor{mygridcolor}{RGB}{0,0,0}
    \definecolor{myhighlightcolor}{RGB}{102,102,102}
    \definecolor{gray1}{RGB}{204,204,204}
    \definecolor{gray2}{RGB}{153,153,153}
\setchessboard{boardfontencoding=LSBC3,
vlabel=false,
hlabel=false,
showmover = false,
    maxfield = h8,
    boardfontsize=30pt,
    boardfontfamily=skaknew,
    pgfstyle=border,
    color=mygridcolor,
    linewidth=0.5pt,
    markboard,
    pgfstyle=color,
    color=mybgcolor,
    backboard,
    color=gray1,
    backregions={h1-h6,h8-h8,f6-f8,d6-e6,a8-b8,a6-a6,e5-e5,c4-d4,f4-f4,h4-h4,a3-c3,e3-g3,b2-c2,e2-f2,a1-b1,g1-g1,d7-d7,g5-g5,c8-d8,c6-c6,d5-d5,f5-f5,a4-a4,d2-d2,g2-g2},
   addblack={qa8,qd7,qc2,qg5},
}
    \makeatletter
    \let\color@endgroupORI\color@endgroup
    \def\color@endgroup{\color@endgroupORI\pgfsetfillopacity{1}}
    \def\cfss@whitefieldmaskcolor{\pgfsetfillopacity{0}\color{white}}
    \def\cfss@blackfieldmaskcolor{\pgfsetfillopacity{0}\color{black}}
    \def\cfss@whitefieldcolor{\pgfsetfillopacity{0}\color{white}}
    \def\cfss@blackfieldcolor{\pgfsetfillopacity{0}\color{black}}
    \makeatother
\scalebox{.5}{
    \chessboard}
\caption{A chessboard with a mod two gradient indicating attacking queens.}
\label{mod2_gradient}
\end{subfigure}
\caption{}
\end{figure}
In the modification of the game, we wish to simplify this illustration as well as try to cover more squares on the board with queens.  Therefore we propose the following additional rule to the two-player game suggested by Noon.  
\begin{quote}
Rule: Let every non-occupied square take on a value given by the number of queens who are directly attacking that square.  If the value is congruent to zero modulo two, then the square is open so a queen may be placed on it, and if the value is congruent to one modulo two the square is closed and we may not place a queen on that square.
\end{quote}
In particular, squares with an even number of attacking queens are open and with an odd number of attacking queens are closed.  Play continues as before until a losing player cannot place another queen on the board.  We will refer to this game as the \textbf{mod 2 $\bm{n}$-queens game}.  Figure~\ref{mod2_gradient} illustrates the same configuration of queen with a now much simpler and more open gradient given by the mod 2 $n$-queens game.  There are many differences between this version of the game, and the traditional $n$-queens problem.  In particular, the order of the placement of the queens now matters.  Further, a maximum of $n^2$ queens can be placed to fill the board, whereas the original game can have at most $n$ queens.
         
In analyzing the game we note there exist three states for any chessboard which we will call complete, unlocked, and locked.   We define these states as follows.
\begin{definition}
A \textbf{complete} chessboard contains $n^2$ queens that have been placed through legal game play.  An \textbf{unlocked} board is one that has less than $n^2$ queens placed on it, and there exist empty open squares in which a queen may be placed.  Finally, a \textbf{locked} board is one that has less than $n^2$ queens, but no legal moves remaining.
\end{definition}

Before we look at some complete chessboards, let us review some chessboard terminology

Each square on a $n\times n$ chessboard will be indexed by an \textbf{ordered pair} $(i,j)$ where $1\leq i \leq n$ indicates the row numbered from top to bottom and $1\leq j \leq n$ indicates the column numbered from left to right.  The \textbf{$\bm{k}$-sum diagonal} is the diagonal running from left to right and bottom to top in which the sum of the indices of each square is $k$ for some integer $1\leq k \leq 2n$.  The \textbf{$\bm{k}$-difference diagonal} is the diagonal running from left to right and top to bottom in which the difference of the indices of each square is $k$ for some integer $-(n-1) \leq k \leq n-1$.  The $(n+1)$-sum diagonal is known as the \textbf{main sum diagonal}, and similarly the $0$-difference diagonal is called the \textbf{main difference diagonal}.

\section{Complete Chessboards}\label{Complete}
By definition, we know that a complete chessboard contains $n^2$ queens, so we ask, is this always possible?

In the case where the size of the board is odd, the answer is yes         
\begin{claim}\label{odd_complete}
If $n$ is an odd positive integer, the $n\times n$ chessboard has a complete solution of $n^2$ queens.
\end{claim}
         
\begin{figure}[ht]
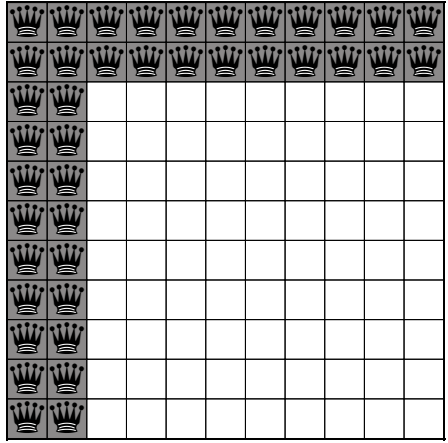

        \centering
    \definecolor{mybgcolor}{RGB}{255,255,255}
    \definecolor{mygridcolor}{RGB}{0,0,0}
    \definecolor{myhighlightcolor}{RGB}{139,137,137}
\setchessboard{boardfontencoding=LSBC3,
vlabel=false,
hlabel=false,
showmover = false,
    maxfield = k11,
    boardfontsize=30pt,
    boardfontfamily=skaknew,
    pgfstyle=border,
    color=mygridcolor,
    linewidth=0.5pt,
    markboard,
    pgfstyle=color,
    color=mybgcolor,
    backboard,
    color=myhighlightcolor,
    backregions={a1-a11,b1-b11,c10-k10,c11-k11},
   addblack={qa1,qa2,qa3,qa4,qa5,qa6,qa7,qa8,qa9,qa10,qa11,
   		qb1,qb2,qb3,qb4,qb5,qb6,qb7,qb8,qb9,qb10,qb11,
		qc10,qd10,qe10,qf10,qg10,qh10,qi10,qj10,qk10,
		qc11,qd11,qe11,qf11,qg11,qh11,qi11,qj11,qk11},
}
    \makeatletter
    \let\color@endgroupORI\color@endgroup
    \def\color@endgroup{\color@endgroupORI\pgfsetfillopacity{1}}
    \def\cfss@whitefieldmaskcolor{\pgfsetfillopacity{0}\color{white}}
    \def\cfss@blackfieldmaskcolor{\pgfsetfillopacity{0}\color{black}}
    \def\cfss@whitefieldcolor{\pgfsetfillopacity{0}\color{white}}
    \def\cfss@blackfieldcolor{\pgfsetfillopacity{0}\color{black}}
    \makeatother
\scalebox{.5}{
    \chessboard
}
       \vspace{-15pt}
        \caption{An unlocked $11\times 11$ chessboard}
        \label{unlocked_11}
\end{figure}

\noindent \textit{Strategy}: Thinking inductively, a $1\times 1$ chessboard may be filled with one queen, and if we want to fill a $3\times 3$ board we need to fill the top two rows and left two columns with queens as well as the $1\times 1$ board in the lower right corner.  More generally, to fill a $n\times n$ board we need to fill the top two rows, the left two columns, and a $(n-2)\times (n-2)$ board in the lower right corner.  Looking at the game position where queens have filled the first two rows and the first two columns of the chessboard as illustrated in Figure~\ref{unlocked_11} in the case $n=11$, we observe that every uncovered square on this chessboard can be attacked by an even number of queens.
Specifically, each square can be attacked vertically, horizontally, and along the difference diagonal by exactly two queens on each line.  Squares to left or on the main diagonal are attacked by four queens, squares on the $n+2$ sum diagonal are attacked by exactly two queens, and all other squares are attacked by zero queens on the sum diagonal.  Therefore game play from the chessboard with queens in this starting configuration is equivalent to the game play on a $(n-2)\times(n-2)$ board with no queens.  Because we can proceed inductively to a $1\times 1$ board, now we only need to show that this starting board can be obtained from a sequence of legal game moves.  

Begin by filling the upper left corner of a $n\times n$ chessboard with the eight queens in the sequence given by Figure~\ref{unlocked_3}.   In the following steps, we will successively fill in the next four squares from the first two rows and column.  We proceed as listed in the steps below and illustrated in Figures~\ref{odd_step1} and~\ref{odd_step2} for an increasing integer $k$ where $1\leq k \leq \frac{n-3}{2}$.
         
        \begin{itemize}
        \item[] Step 1: Place a queen in square $(2,2k+2)$ and then in square $(1, 2k+2)$.
        \item[] Step 2: Place a queen in square $(2k+2, 1)$ and then in square $(2k+2,2)$.
        \item[] Step 3: Place a queen in square $(1,2k+3)$ and then in square $(2, 2k+3)$.
        \item[] Step 4: Place a queen in square $(2k+3, 2)$ and then in square $(2k+3,1)$.
        \end{itemize}

\begin{figure}[ht]
        \centering
        \begin{subfigure}{.3\textwidth}
         \definecolor{mybgcolor}{RGB}{255,255,255}
    \definecolor{mygridcolor}{RGB}{0,0,0}
    \definecolor{myhighlightcolor}{RGB}{139,137,137}
\setchessboard{boardfontencoding=LSBC3,
vlabel=false,
hlabel=false,
showmover = false,
    maxfield = g7,
    boardfontsize=30pt,
    boardfontfamily=skaknew,
    pgfstyle=border,
    color=mygridcolor,
    linewidth=0.5pt,
    markboard,
    pgfstyle=color,
    color=mybgcolor,
    backboard,
    color=myhighlightcolor,
    backregions={a1-a4, b1-b3, c3-c3, d7-d7,d2-d2, e5-e7, e1-e1, f6-f7, f4-f4, g6-g7, g3-g3},
}
    \makeatletter
    \let\color@endgroupORI\color@endgroup
    \def\color@endgroup{\color@endgroupORI\pgfsetfillopacity{1}}
    \def\cfss@whitefieldmaskcolor{\pgfsetfillopacity{0}\color{white}}
    \def\cfss@blackfieldmaskcolor{\pgfsetfillopacity{0}\color{black}}
    \def\cfss@whitefieldcolor{\pgfsetfillopacity{0}\color{white}}
    \def\cfss@blackfieldcolor{\pgfsetfillopacity{0}\color{black}}
    \makeatother
\scalebox{.5}{
    \chessboard[pgfstyle=text,
         text=1,
         markregions={a7-a7},
         text=2,
         markregions={b5-b5},
         text=3,
         markregions={b6-b6},
         text=4,
         markregions={c6-c6},
         text=5,
         markregions={a6-a6},
         text=6,
         markregions={a5-a5},
         text=7,
         markregions={c7-c7},
         text=8,
         markregions={b7-b7}]
}        
%
        \vspace{-15pt}
        \caption{Place the first eight queens.}
        \label{unlocked_3}
        \end{subfigure}
        \qquad
         
        \begin{subfigure}{.4\textwidth}
        \centering
                 \definecolor{mybgcolor}{RGB}{255,255,255}
    \definecolor{mygridcolor}{RGB}{0,0,0}
    \definecolor{myhighlightcolor}{RGB}{139,137,137}
\setchessboard{boardfontencoding=LSBC3,
vlabel=false,
hlabel=false,
showmover = false,
    maxfield = g7,
    boardfontsize=30pt,
    boardfontfamily=skaknew,
    pgfstyle=border,
    color=mygridcolor,
    linewidth=0.5pt,
    markboard,
    pgfstyle=color,
    color=mybgcolor,
    backboard,
    color=myhighlightcolor,
    backregions={a1-a2, b1-b1, c1-c1, f7-f7, g5-g7, a5-a7, b5-b7, c7-c7, c6-c6},
    addblack={qa7, qb5, qb6, qc6, qa6, qa5, qc7, qb7},
}
    \makeatletter
    \let\color@endgroupORI\color@endgroup
    \def\color@endgroup{\color@endgroupORI\pgfsetfillopacity{1}}
    \def\cfss@whitefieldmaskcolor{\pgfsetfillopacity{0}\color{white}}
    \def\cfss@blackfieldmaskcolor{\pgfsetfillopacity{0}\color{black}}
    \def\cfss@whitefieldcolor{\pgfsetfillopacity{0}\color{white}}
    \def\cfss@blackfieldcolor{\pgfsetfillopacity{0}\color{black}}
    \makeatother
\scalebox{.5}{
    \chessboard[pgfstyle=text,
         text=1,
         markregions={d6-d6},
         text=2,
         markregions={d7-d7},
         text=3,
         markregions={a4-a4},
         text=4,
         markregions={b4-b4},
         text=5,
         markregions={e7-e7},
         text=6,
         markregions={e6-e6},
         text=7,
         markregions={b3-b3},
         text=8,
         markregions={a3-a3}]
}  
        \vspace{-15pt}
        \caption{Place the next eight queens.}
        \label{odd_step1}
        \end{subfigure}
        \qquad
        \begin{subfigure}{.4\textwidth}
        \centering
                   \definecolor{mybgcolor}{RGB}{255,255,255}
    \definecolor{mygridcolor}{RGB}{0,0,0}
    \definecolor{myhighlightcolor}{RGB}{139,137,137}
\setchessboard{boardfontencoding=LSBC3,
vlabel=false,
hlabel=false,
showmover = false,
    maxfield = g7,
    boardfontsize=30pt,
    boardfontfamily=skaknew,
    pgfstyle=border,
    color=mygridcolor,
    linewidth=0.5pt,
    markboard,
    pgfstyle=color,
    color=mybgcolor,
    backboard,
    color=myhighlightcolor,
    backregions={a3-a7, b3-b7, c7-e7, c6-e6},
    addblack={qa7, qb5, qb6, qc6, qa6, qa5, qc7, qb7, qd6, qd7, qb4, qa4, qe7, qe6, qb3, qa3},
}
    \makeatletter
    \let\color@endgroupORI\color@endgroup
    \def\color@endgroup{\color@endgroupORI\pgfsetfillopacity{1}}
    \def\cfss@whitefieldmaskcolor{\pgfsetfillopacity{0}\color{white}}
    \def\cfss@blackfieldmaskcolor{\pgfsetfillopacity{0}\color{black}}
    \def\cfss@whitefieldcolor{\pgfsetfillopacity{0}\color{white}}
    \def\cfss@blackfieldcolor{\pgfsetfillopacity{0}\color{black}}
    \makeatother
\scalebox{.5}{
    \chessboard[pgfstyle=text,
         text=1,
         markregions={f6-f6},
         text=2,
         markregions={f7-f7},
         text=3,
         markregions={a2-a2},
         text=4,
         markregions={b2-b2},
         text=5,
         markregions={g7-g7},
         text=6,
         markregions={g6-g6},
         text=7,
         markregions={b1-b1},
         text=8,
         markregions={a1-a1}]
}       
        \vspace{-15pt}
        \caption{Place the final eight queens.}
        \label{odd_step2}
        \end{subfigure}
        \caption{Steps in the construction of a complete board for $n$ odd.}
        \label{XX}
        \end{figure}
        
Once this process is finished, we repeat on the $(n-2)\times (n-2)$ chessboard until we have reduced to the one empty square in the lower right corner which can then be filled with a queen.  Not all chessboards, however, will have complete solutions.         
         
\begin{claim}\label{even_not_complete}
If $n$ is an even positive integer, the $n\times n$ chessboard does not have a complete solution of $n^2$ queens.
\end{claim}

\noindent \textit{Contradiction:}  To see how this is true, suppose to the contrary, that there is an open square in which to place $n^2$th queen on a $n\times n$ chessboard, and further suppose that square is on the outer edge of the chessboard.  Each square on the edge has $n-1$ queens attacking along the horizontal line and $n-1$ queens attacking along the vertical line.  If we chose a corner square, it only has one diagonal containing $n-1$ attacking queens.  As you move from a corner square along the edge, the larger diagonal decreases by one and the lesser diagonal increases by one which always keeps the sum of the diagonals of a square on the edge equal to $n-1$.  Thus, the open square cannot be one of the squares on the edge of the chessboard because the number of attacking queens is $(n - 1) + (n - 1) + (n - 1) \equiv_2 1$.  We know the outer ring is filled with queens, so we remove this ring without changing the parity of the chessboard, as each square in the middle is attacked by exactly eight queens from the outer ring.  We are left with $(n-2)\times(n-2)$ chessboard.  As $n-2$ is still even, we repeat, removing the outer rings, until we are left with the contradiction, a $2\times 2$ chessboard which cannot be filled with four queens.

Although we cannot have a complete solution in the case of an even length chessboard, we can identify a locked solution that comes close.  

\begin{claim}\label{even_complete}
If $n$ is an even positive integer, $n^2-2$ queens can be placed legally on a $n\times n$ chessboard.
\end{claim}

One such construction could follow similarly to the odd case by applying induction and filling the topmost two rows and leftmost two columns, so we leave it as a challenge.  This case is interesting, however, and we pose the following question.

\begin{ques}
Is the solution of $n^2-2$ queens on an $n\times n$ chessboard for $n$ an even positive integer a maximal locked solution?
\end{ques}

Empirical data and a computer simulation for $n=4$ suggest that the answer is yes.  However, even though the program could run the simulation for $n=4$ in eleven minutes, it would take approximately 22 years to run for $n=6$, so verifying computationally is not an option.

In the next section, we consider locked chessboards.         
         
\section{Locked Chessboards}\label{Locked}
         
Another interesting game state is a locked chessboard.  For strategy reasons, a player would be interested to know of locked solutions in order to lead an opponent into such a solution or avoid them himself.  For small $n=1, 2$ or $3$, there is a simple first-player win strategy, that is, placing only one queen as seen in Figures~\ref{locked_2} and \ref{locked_3}.
\begin{figure}[H]
\centering
\begin{subfigure}[t]{.46\textwidth}
\centering
\vspace{-10pt}
    \definecolor{mybgcolor}{RGB}{255,255,255}
    \definecolor{mygridcolor}{RGB}{0,0,0}
    \definecolor{myhighlightcolor}{RGB}{139,137,137}
\setchessboard{boardfontencoding=LSBC3,
vlabel=false,
hlabel=false,
showmover = false,
    maxfield = b2,
    boardfontsize=30pt,
    boardfontfamily=skaknew,
    pgfstyle=border,
    color=mygridcolor,
    linewidth=0.5pt,
    markboard,
    pgfstyle=color,
    color=mybgcolor,
    backboard,
    color=myhighlightcolor,
    backregions={a1-b2},
   addblack={qa2},
}
    \makeatletter
    \let\color@endgroupORI\color@endgroup
    \def\color@endgroup{\color@endgroupORI\pgfsetfillopacity{1}}
    \def\cfss@whitefieldmaskcolor{\pgfsetfillopacity{0}\color{white}}
    \def\cfss@blackfieldmaskcolor{\pgfsetfillopacity{0}\color{black}}
    \def\cfss@whitefieldcolor{\pgfsetfillopacity{0}\color{white}}
    \def\cfss@blackfieldcolor{\pgfsetfillopacity{0}\color{black}}
    \makeatother
\scalebox{.5}{
    \chessboard
}
\vspace{-15pt}
\caption{A locked $2\times 2$ chessboard}
\label{locked_2}
\end{subfigure}
\qquad
\begin{subfigure}[t]{.46\textwidth}
\centering
\vspace{-10pt}
    \definecolor{mybgcolor}{RGB}{255,255,255}
    \definecolor{mygridcolor}{RGB}{0,0,0}
    \definecolor{myhighlightcolor}{RGB}{139,137,137}
\setchessboard{boardfontencoding=LSBC3,
vlabel=false,
hlabel=false,
showmover = false,
    maxfield = c3,
    boardfontsize=30pt,
    boardfontfamily=skaknew,
    pgfstyle=border,
    color=mygridcolor,
    linewidth=0.5pt,
    markboard,
    pgfstyle=color,
    color=mybgcolor,
    backboard,
    color=myhighlightcolor,
    backregions={a1-c3},
   addblack={qb2},
}
    \makeatletter
    \let\color@endgroupORI\color@endgroup
    \def\color@endgroup{\color@endgroupORI\pgfsetfillopacity{1}}
    \def\cfss@whitefieldmaskcolor{\pgfsetfillopacity{0}\color{white}}
    \def\cfss@blackfieldmaskcolor{\pgfsetfillopacity{0}\color{black}}
    \def\cfss@whitefieldcolor{\pgfsetfillopacity{0}\color{white}}
    \def\cfss@blackfieldcolor{\pgfsetfillopacity{0}\color{black}}
    \makeatother
\scalebox{.5}{
    \chessboard
}
\vspace{-15pt}
\caption{A locked $3\times 3$ chessboard}
\label{locked_3}
\end{subfigure}
\end{figure}

Of course, these small examples cannot be generalized and locked boards with only one queen do not exist for $n>3$.  Let's look at a more complex locked chessboard.  Again, as we did for complete boards, we separate the chessboards into two classes by their even or odd lengths.  

\begin{figure}[H]
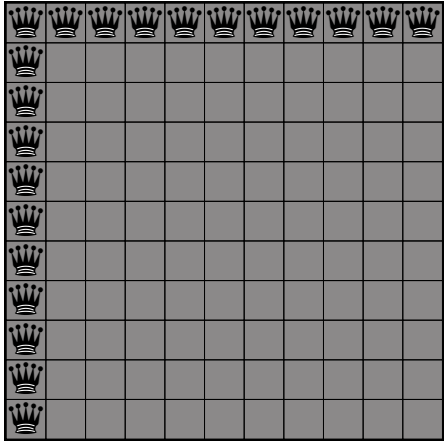
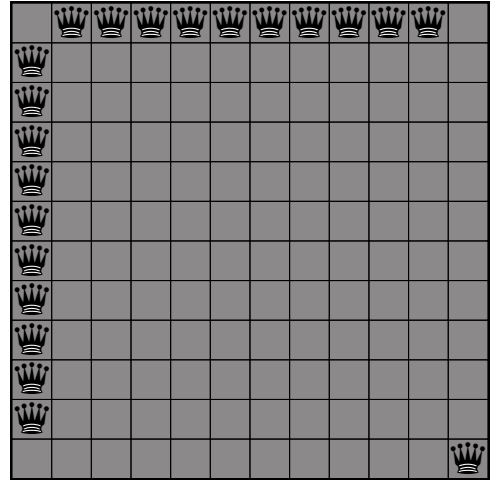

\centering
\begin{subfigure}[t]{.46\textwidth}
\centering
\vspace{-10pt}
    \definecolor{mybgcolor}{RGB}{255,255,255}
    \definecolor{mygridcolor}{RGB}{0,0,0}
    \definecolor{myhighlightcolor}{RGB}{139,137,137}
\setchessboard{boardfontencoding=LSBC3,
vlabel=false,
hlabel=false,
showmover = false,
    maxfield = k11,
    boardfontsize=30pt,
    boardfontfamily=skaknew,
    pgfstyle=border,
    color=mygridcolor,
    linewidth=0.5pt,
    markboard,
    pgfstyle=color,
    color=mybgcolor,
    backboard,
    color=myhighlightcolor,
    backregions={a1-k11},
   addblack={qa1,qa2,qa3,qa4,qa5,qa6,qa7,qa8,qa9,qa10,qa11
   		,qb11,qc11,qd11,qe11,qf11,qg11,qh11,qi11,qj11,qk11},
}
    \makeatletter
    \let\color@endgroupORI\color@endgroup
    \def\color@endgroup{\color@endgroupORI\pgfsetfillopacity{1}}
    \def\cfss@whitefieldmaskcolor{\pgfsetfillopacity{0}\color{white}}
    \def\cfss@blackfieldmaskcolor{\pgfsetfillopacity{0}\color{black}}
    \def\cfss@whitefieldcolor{\pgfsetfillopacity{0}\color{white}}
    \def\cfss@blackfieldcolor{\pgfsetfillopacity{0}\color{black}}
    \makeatother
\scalebox{.5}{
    \chessboard
}
\vspace{-15pt}
\caption{A locked $11\times 11$ chessboard}
\label{locked_11}
\end{subfigure}
\qquad
\begin{subfigure}[t]{.46\textwidth}
\centering
\vspace{-10pt}
 \definecolor{mybgcolor}{RGB}{255,255,255}
    \definecolor{mygridcolor}{RGB}{0,0,0}
    \definecolor{myhighlightcolor}{RGB}{139,137,137}
\setchessboard{boardfontencoding=LSBC3,
vlabel=false,
hlabel=false,
showmover = false,
    maxfield = l12,
    boardfontsize=30pt,
    boardfontfamily=skaknew,
    pgfstyle=border,
    color=mygridcolor,
    linewidth=0.5pt,
    markboard,
    pgfstyle=color,
    color=mybgcolor,
    backboard,
    color=myhighlightcolor,
    backregions={a1-l12},
   addblack={ql1,qa2,qa3,qa4,qa5,qa6,qa7,qa8,qa9,qa10,qa11
   		,qb12,qc12,qd12,qe12,qf12,qg12,qh12,qi12,qj12,qk12},
}
    \makeatletter
    \let\color@endgroupORI\color@endgroup
    \def\color@endgroup{\color@endgroupORI\pgfsetfillopacity{1}}
    \def\cfss@whitefieldmaskcolor{\pgfsetfillopacity{0}\color{white}}
    \def\cfss@blackfieldmaskcolor{\pgfsetfillopacity{0}\color{black}}
    \def\cfss@whitefieldcolor{\pgfsetfillopacity{0}\color{white}}
    \def\cfss@blackfieldcolor{\pgfsetfillopacity{0}\color{black}}
    \makeatother
\scalebox{.5}{ 
    \chessboard
}
\vspace{-15pt}
\caption{A locked $12\times 12$ chessboard}
\label{locked_12}
\end{subfigure}
\caption{Examples of even-sized and odd-sized locked boards}
\label{general_locked}
\end{figure}

Consider the chessboards where queens fill the top row and leftmost column for $n$ odd or almost fill the top row and leftmost column for $n$ even.  These game positions are illustrated in Figure~\ref{general_locked}.
        
\begin{claim}\label{odd_locked}
For $n$ an odd positive integer, the chessboard with queens exactly on squares in the set $\{(1,i), (i,1) | 1\leq i \leq n\}$ is a legal, locked game position.
\end{claim}

\noindent \textit{Strategy:} We need to check that this board is locked as well as check that the queens can be placed in this configuration by a legal set of game moves.  First we confirm that an odd number of queens attack every uncovered square.  Observe that every unoccupied square has exactly one queen which can attack vertically and exactly one queen that can attack horizontally.  On the positive diagonal, if the uncovered square is on or to the left of the main sum diagonal exactly two queens can attack along the positive diagonal.  However, if the uncovered square is the right of the main diagonal, no queens can attack along the positive diagonal.  Further every empty square can be attacked by exactly one queen on the difference diagonal, so all uncovered square are attacked by either exactly five or exactly three queens and hence the board is locked.

\begin{figure}[ht]
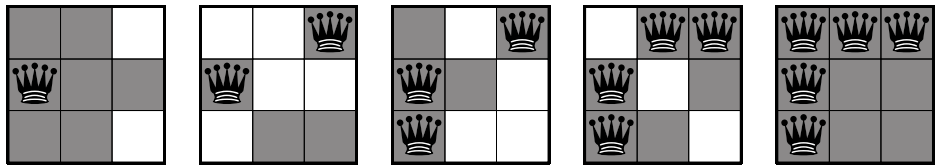

\centering
\hspace{-58pt}
\begin{subfigure}[t]{.11\textwidth}
    \definecolor{mybgcolor}{RGB}{255,255,255}
    \definecolor{mygridcolor}{RGB}{0,0,0}
    \definecolor{myhighlightcolor}{RGB}{139,137,137}
\setchessboard{boardfontencoding=LSBC3,
vlabel=false,
hlabel=false,
showmover = false,
    maxfield = c3,
    boardfontsize=30pt,
    boardfontfamily=skaknew,
    pgfstyle=border,
    color=mygridcolor,
    linewidth=0.5pt,
    markboard,
    pgfstyle=color,
    color=mybgcolor,
    backboard,
    color=myhighlightcolor,
    backregions={a1-a3,b2-c2,b1-b3},
   addblack={qa2},
}
    \makeatletter
    \let\color@endgroupORI\color@endgroup
    \def\color@endgroup{\color@endgroupORI\pgfsetfillopacity{1}}
    \def\cfss@whitefieldmaskcolor{\pgfsetfillopacity{0}\color{white}}
    \def\cfss@blackfieldmaskcolor{\pgfsetfillopacity{0}\color{black}}
    \def\cfss@whitefieldcolor{\pgfsetfillopacity{0}\color{white}}
    \def\cfss@blackfieldcolor{\pgfsetfillopacity{0}\color{black}}
    \makeatother
\scalebox{.65}{
    \chessboard
}
\vspace{-15pt}
\end{subfigure}
\qquad
\begin{subfigure}[t]{.11\textwidth}
  \definecolor{mybgcolor}{RGB}{255,255,255}
    \definecolor{mygridcolor}{RGB}{0,0,0}
    \definecolor{myhighlightcolor}{RGB}{139,137,137}
\setchessboard{boardfontencoding=LSBC3,
vlabel=false,
hlabel=false,
showmover = false,
    maxfield = c3,
    boardfontsize=30pt,
    boardfontfamily=skaknew,
    pgfstyle=border,
    color=mygridcolor,
    linewidth=0.5pt,
    markboard,
    pgfstyle=color,
    color=mybgcolor,
    backboard,
    color=myhighlightcolor,
    backregions={a2-a2,c3-c3,b1-c1},
   addblack={qa2,qc3},
}
    \makeatletter
    \let\color@endgroupORI\color@endgroup
    \def\color@endgroup{\color@endgroupORI\pgfsetfillopacity{1}}
    \def\cfss@whitefieldmaskcolor{\pgfsetfillopacity{0}\color{white}}
    \def\cfss@blackfieldmaskcolor{\pgfsetfillopacity{0}\color{black}}
    \def\cfss@whitefieldcolor{\pgfsetfillopacity{0}\color{white}}
    \def\cfss@blackfieldcolor{\pgfsetfillopacity{0}\color{black}}
    \makeatother
\scalebox{.65}{
    \chessboard
}
\vspace{-15pt}
\end{subfigure}
\qquad
\begin{subfigure}[t]{.11\textwidth}
 \definecolor{mybgcolor}{RGB}{255,255,255}
    \definecolor{mygridcolor}{RGB}{0,0,0}
    \definecolor{myhighlightcolor}{RGB}{139,137,137}
\setchessboard{boardfontencoding=LSBC3,
vlabel=false,
hlabel=false,
showmover = false,
    maxfield = c3,
    boardfontsize=30pt,
    boardfontfamily=skaknew,
    pgfstyle=border,
    color=mygridcolor,
    linewidth=0.5pt,
    markboard,
    pgfstyle=color,
    color=mybgcolor,
    backboard,
    color=myhighlightcolor,
    backregions={a1-a3,b2-b2,c3-c3},
   addblack={qa2,qc3,qa1},
}
    \makeatletter
    \let\color@endgroupORI\color@endgroup
    \def\color@endgroup{\color@endgroupORI\pgfsetfillopacity{1}}
    \def\cfss@whitefieldmaskcolor{\pgfsetfillopacity{0}\color{white}}
    \def\cfss@blackfieldmaskcolor{\pgfsetfillopacity{0}\color{black}}
    \def\cfss@whitefieldcolor{\pgfsetfillopacity{0}\color{white}}
    \def\cfss@blackfieldcolor{\pgfsetfillopacity{0}\color{black}}
    \makeatother
\scalebox{.65}{
    \chessboard
}
\vspace{-15pt}
\end{subfigure}
\qquad
\begin{subfigure}[t]{.11\textwidth}
  \definecolor{mybgcolor}{RGB}{255,255,255}
    \definecolor{mygridcolor}{RGB}{0,0,0}
    \definecolor{myhighlightcolor}{RGB}{139,137,137}
\setchessboard{boardfontencoding=LSBC3,
vlabel=false,
hlabel=false,
showmover = false,
    maxfield = c3,
    boardfontsize=30pt,
    boardfontfamily=skaknew,
    pgfstyle=border,
    color=mygridcolor,
    linewidth=0.5pt,
    markboard,
    pgfstyle=color,
    color=mybgcolor,
    backboard,
    color=myhighlightcolor,
    backregions={a1-a2,b3-c3,c2-c2,b1-b1},
   addblack={qa2,qc3,qa1,qb3},
}
    \makeatletter
    \let\color@endgroupORI\color@endgroup
    \def\color@endgroup{\color@endgroupORI\pgfsetfillopacity{1}}
    \def\cfss@whitefieldmaskcolor{\pgfsetfillopacity{0}\color{white}}
    \def\cfss@blackfieldmaskcolor{\pgfsetfillopacity{0}\color{black}}
    \def\cfss@whitefieldcolor{\pgfsetfillopacity{0}\color{white}}
    \def\cfss@blackfieldcolor{\pgfsetfillopacity{0}\color{black}}
    \makeatother
\scalebox{.65}{
    \chessboard
}
\vspace{-15pt}
\end{subfigure}
\qquad
\begin{subfigure}[t]{.11\textwidth}
 \definecolor{mybgcolor}{RGB}{255,255,255}
    \definecolor{mygridcolor}{RGB}{0,0,0}
    \definecolor{myhighlightcolor}{RGB}{139,137,137}
\setchessboard{boardfontencoding=LSBC3,
vlabel=false,
hlabel=false,
showmover = false,
    maxfield = c3,
    boardfontsize=30pt,
    boardfontfamily=skaknew,
    pgfstyle=border,
    color=mygridcolor,
    linewidth=0.5pt,
    markboard,
    pgfstyle=color,
    color=mybgcolor,
    backboard,
    color=myhighlightcolor,
    backregions={a1-c3},
   addblack={qa2,qc3,qa1,qb3, qa3},
}
    \makeatletter
    \let\color@endgroupORI\color@endgroup
    \def\color@endgroup{\color@endgroupORI\pgfsetfillopacity{1}}
    \def\cfss@whitefieldmaskcolor{\pgfsetfillopacity{0}\color{white}}
    \def\cfss@blackfieldmaskcolor{\pgfsetfillopacity{0}\color{black}}
    \def\cfss@whitefieldcolor{\pgfsetfillopacity{0}\color{white}}
    \def\cfss@blackfieldcolor{\pgfsetfillopacity{0}\color{black}}
    \makeatother
\scalebox{.65}{
    \chessboard
}
\vspace{-15pt}
\end{subfigure}
\caption{A sequence of game positions to create a locked $3\times 3$ chessboard}
\label{locked_n=3}
\end{figure}

Next, we need to show that this game board can arise from a sequence of legal placements of queens onto the board.  First consider the case where $n=3$.  Figure~\ref{locked_n=3} illustrates a sequence of five moves to lock a $3\times 3$ chessboard by playing queens on the first row and column.

We want to extend the game play in Figure~\ref{locked_n=3} to occur in a chessboard of any size.  Our strategy will be to repeat the first four moves of  Figure~\ref{locked_n=3} for the next two vacant squares in the first row and the corresponding two vacant squares in the first column.  Suppose we wish to fill the four squares from the set $\{(1,i), (1,i+1), (i,1), (i+1,1)\}$ with queens.  In each stage will fill two squares in column one and two squares in row one, so the total number of queens attacking either horizontally or vertically must be even and further there are no queens which can attack from either of the two diagonals.  Squares $(i,1)$ and $(1, i+1)$ are not on a line vertically, horizontally, or diagonally, so a queen can be placed on each of these squares.  The remaining squares $(i+1,1)$ and $(1,i)$ can be attacked by both queens that were just played and also are not on a line with each other vertically, horizontally, or diagonally.  We can fill these squares with queens.  After repeating this process $\frac{n-1}{2}$ times, each square in the first row and first column, except the square in the upper left corner, contains a queen.  As the square in the upper left can be attacked by all $2(n-1)$ queens on the board, we can place the final queen in this position and hence lock the board.

We can find a locked board in the case where $n$ is even as well.

\begin{claim}
For $n$ an even positive integer and $n>2$, the chessboard with queens on squares in the set $\{(1,i), (i,1) | 2\leq i \leq n-1\} \cup \{(n,n)\}$ is a legal, locked game position.
\end{claim}

We encourage you to check that this chessboard is locked and arises from a sequence of legal game moves.

We can now make the following observation.
\begin{corollary}
When $n$ is an odd positive integer, at most $2n -1$ queens are needed to lock the board, and
when $n$ is an even positive integer, at most $2n-3$ queens are needed to lock the board.
\end{corollary}

This corollary gives an upper bound on the number of queens need to lock an even or odd chessboard, but is this bound strict?  We know for very small values of $n$ the chessboard can be locked with fewer queens, but what about larger values?  We pose the next question.

\begin{ques}
What is the minimum number of queens needed to lock a $n\times n$ chessboard?
\end{ques}

Next, we will look at different versions of the game with alternate rules and the connections between them.         
         
\section{In the Alternate Universe}\label{Alternate}
Thus far, we have seen that there is a difference in outcomes for chessboards with even and odd length, that is, with parity zero or one modulo two respectively.  The rules of the game as well could be modified along the same lines.  Let's look at the alternative version of the game where queens can be placed on squares that are attacked by an \textit{odd} number of queens.  We will call such a game an \textbf{alternate universe mod 2 $\bm{n}$-queens game}.
         
        The first observation of course is that generally this cannot happen as an empty board has no squares attacked by on odd number of queens, so we will modify the rules again to start with one queen already on the board before we commence play.  From a strategy standpoint, in this alternate universe game a new queen cannot be placed unless she is ``covered" by an odd number of queens.
         
        A second way to play the game would be to start with $n^2$ queens on the board and try to remove the queens one by one following the rule that a queen may be removed as long as she can be attacked by an even number of queens.  We call such a game an \textbf{complementary mod 2 $\bm{n}$-queens game}.  The relationship between alternate universe, complementary games, and standard mod 2 $n$-queens games is dependent on the size of the chessboard.
 \begin{claim}
 When $n$ is an even positive integer, the alternate universe mod 2 $n$-queens games are in bijection with complementary mod 2 $n$-queens games.
\end{claim}
         
        Earlier we saw that if an even length chessboard is fully covered with queens, each square is attacked by an odd number of queens.  The initial step for both games in the same.  We place one queen in the alternate universe and remove one queen from the same square in the full complementary game.  While playing the complementary game, queens can only be removed if the number of queens attacking is even and hence because of the parity of the board, it means the number of empty squares ``attacking" must be odd.  By taking the complement of the board, removing queens to make empty squares, and placing queens on otherwise empty squares, we can go back and forth from the complementary game to the alternate universe game.
         
        It is interesting to note that the alternate universe game on the odd length chessboards does not have the same relationship.  In fact, we have the following:
         
\begin{claim}
When $n$ is an odd positive integer, the complementary mod 2 $n$-queens games are in bijection with standard mod 2 $n$-queens games.
\end{claim}
         
        In this case, the beginning parity for every square on a complementary game board is even. We obtain the bijection by placing a queen in the same square of the standard game from which we removed it in the complementary game.  

Thus for each case, even or odd length chessboard, there is really only two different games to be played, standard and alternate universe, as the complementary game can be completely reconstructed from one of the other two games.

We conclude this section with an open-ended question.

\begin{ques}
What other $n$-queens games can be played on a $n\times n$ chessboard?  In particular, how could we adapt play for a \textbf{mod k $\bm{n}$-queens game}?
\end{ques}
        
In the next section we will discuss the complexity of the game and give examples of a game tree and a graph created from this tree.

\section{Game Trees, Graphs, and Complexity}\label{Graph}

\begin{figure}[ht]
\hspace{-18pt}
\includegraphics[scale=.58]{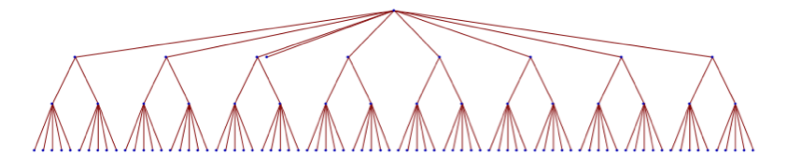}
\caption{Partial game tree for $n=3$}
\label{game_tree}
\end{figure}

Chess and other games played with chess pieces are complex games.  We will discuss several standard measures of complexity of a game for the mod 2 $n$-queens game.  A \textbf{game tree} is the rooted tree whose root is the empty board and whose leaves are locked or complete boards such that the directed paths from the root to a leaf display all possible games.  Figure~\ref{game_tree} shows the first three levels of the game tree in the case where $n=3$.
Because there are $n^2$ squares on the board, we can get an upper bound for the number of leaf nodes or \textbf{game tree size} at $n^2!$.

\begin{figure}[ht]
\hspace{20pt}
\includegraphics[scale=.6]{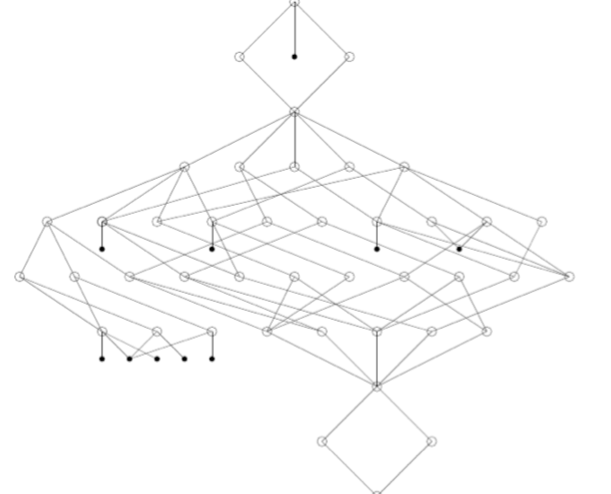}
\caption{Game graph for $n=3$}
\label{fullgraph}
\end{figure}

Many of the game boards represented by these nodes are identical or symmetric by rotation or reflection.  In fact, any given game board has up to eight symmetric boards. To simplify the structure of this tree, we propose to combine nodes describing boards which are equivalent by symmetry.  In this case, we have fewer nodes, but the tree has become a graph.  Figure~\ref{fullgraph} shows the full graph for $n = 3$.

        This graph is still quite large even for smallest non-trivial value for $n$.  Observe that that the first node is the empty board, the last board is the full board, and each of the solid nodes correspond to a locked board.

We are also interested in the complexity of the game.  Because, for some $n$, all possible squares may be filled, we have a finite set of $n^2$ squares; any subset of which could be filled with queens.  This set of all subsets of $n^2$ elements is the power set which has size of $2^{n^2}$.  Since defining rules for the board limits the number of game boards, $2^{n^2}$ is the upper bound for the \textbf{state-space complexity} or number of legal game positions.  Therefore we have a complexity $O(2^{n^2})$ and this game is in EXPTIME.

 We conclude with one final question.
 \begin{ques}
 Clearly, the upper bounds for the game tree size and the state-space complexity are not firm for all $n$.  Can these be improved?
 \end{ques}

\newcommand{\journal}[6]{{\sc #1,} #2, {\it #3} {\bf #4} (#5), #6.}
\newcommand{\book}[4]{{\sc #1,} ``#2,'' #3, #4.}
\newcommand{\bookf}[5]{{\sc #1,} ``#2,'' #3, #4, #5.}
\newcommand{\thesis}[4]{{\sc #1,} ``#2,'' Doctoral dissertation, #3, #4.}
\newcommand{\masters}[4]{{\sc #1,} ``#2,'' Thesis, #3, #4.}
\newcommand{\springer}[4]{{\sc #1,} ``#2,'' Lecture Notes in Math., 
                          Vol.\ #3, Springer-Verlag, Berlin, #4.}
\newcommand{\preprint}[3]{{\sc #1,} #2, preprint #3.}
\newcommand{\progress}[2]{{\sc #1,} #2, work in progress.}
\newcommand{\archive}[3]{{\sc #1,} #2, {\bf #3}.}
\newcommand{\unpublished}[1]{{\sc #1,} unpublished.}
\newcommand{\unpublisheddate}[2]{{\sc #1,} unpublished #2.}
\newcommand{\preparation}[2]{{\sc #1,} #2, in preparation.}
\newcommand{\appear}[3]{{\sc #1,} #2, to appear in {\it #3}.}
\newcommand{\submitted}[4]{{\sc #1,} #2, submitted to {\it #3}, #4.}
\newcommand{\AdvancesinMathematics}{Adv.\ Math.}
\newcommand{\DiscreteComputationalGeometry}{Discrete Comput.\ Geom.}
\newcommand{\DiscreteMath}{Discrete Math.}
\newcommand{\EuropeanJournalofCombinatorics}{European J.\ Combin.}
\newcommand{\JCTA}{J.\ Combin.\ Theory Ser.\ A}
\newcommand{\JCTB}{J.\ Combin.\ Theory Ser.\ B}
\newcommand{\JournalofAlgebraicCombinatorics}{J.\ Algebraic Combin.}
\newcommand{\communication}[1]{{\sc #1,} personal communication.}
\newcommand{\website}[2]{{\sc #1}, #2.}

\newcommand{\collection}[9]{{\sc #1,} #2, 
           in {\it #3} (#4), {\it #5},
           #6, #7, #8, #9.}

 \end{document}